\documentclass{amsart}
\usepackage{a4,amssymb}
\def\bar{\overline}
\newcommand{\bbas}{\begin{eqnarray*}}
\newcommand{\eeas}{\end{eqnarray*}}
\newcommand{\bbar}{\begin{array}}
\newcommand{\eear}{\end{array}}
\newcommand{\bbs}{\begin{displaymath}}
\newcommand{\ees}{\end{displaymath}}
\newcommand{\bb}{\begin{equation}}
\newcommand{\eqbb}{\begin{equation}}
\def\ee{\end{equation}}
\def\eqee{\end{equation}}
\def\eea{\end{eqnarray}}
\def\bba{\begin{eqnarray}}

\newtheorem{thm}{Theorem}

\def\lora{\longrightarrow}
\def\ot{\otimes}

\def\loma{\longmapsto}
\def\proof{{\it Proof.\ }}

\def\eee{\rule{.75ex}{1.5ex}\vskip1ex}

\def\Ker{\mbox{\rm Ker}{}}
\def\Im{\mbox{\rm Im}{}}

\def\Hom{\mbox{\rm Hom{}}}
\def\End{\mbox{\rm End{}}}

\def\fun#1{{\mathcal #1}}

\newcommand{\rk}{\mbox{\rm rank}}

\renewcommand{\dim}{\mbox{\rm dim}}

\font\Fraktur=eufm10 scaled\magstep1          
   \newcommand{\fraktur}[1]{\mbox{\Fraktur #1}}  %
   \font\Fraktu=eufm7 scaled\magstep1            
   \newcommand{\fraktu}[1]{\mbox{\Fraktu #1}}    %
   \font\Frakt=eufm5 scaled\magstep1             
  \newcommand{\frakt}[1]{\mbox{\Frakt #1}}      %
   \def\fr#1{\mathchoice{\fraktur {#1}}            
                        {\fraktur {#1}}            
                        {\fraktu {#1}}             
                        {\frakt {#1}}  }           
\newcommand{\Ss}{\fr S}

\def\H{{\mathcal H}}

\def\db{{\mathchoice{\mbox{\rm db}}
                    {\mbox{\rm db}}
                    {\mbox{\scriptsize\rm db}}
                    {\mbox{\tiny\rm db}} }}

\def\ev{{\mathchoice{\mbox{\rm ev}}
                    {\mbox{\rm ev}}
                    {\mbox{\scriptsize\rm ev}}
                    {\mbox{\tiny\rm ev}} }}

\def\id{{\mathchoice{\mbox{\rm id}}
                    {\mbox{\rm id}}
                    {\mbox{\scriptsize\rm id}}
                    {\mbox{\tiny\rm id}} }}

\def\det{\mbox{\rm det{}}}

\def\sdet{\mbox{\rm Ber$_q$}}
\def\comod{\mbox{-comod}}

\def\myperaddress{Institute of Mathematics, P.O.Box 631, 10000 Boho, Hanoi}
\newtheorem{fact}{}[section]

\def\k{k}
\begin{document}
\title{The homological determinant of quatum groups of type $A$}
\author{Ph\`ung H$\grave{\mbox{\^ o}}$ Hai}
\address{\myperaddress}
\thanks{This work is supported by the National Program of Basic Sciences Research of Vietnam.}
\maketitle

\bibliographystyle{plain}

\section*{Introduction}
Let $V$ be a vector space over a field $k$ and $G=GL(V)$, the general linear group.
Let $n=\dim_kV$.  It is well known 
that elements of $G$ acts on the $n$-th homogeneous component of the exterior  over $V$ by
means of the determinant. More precisely, let $x_1,x_2,\ldots, x_n$ be a basis of $V$. Then $\Lambda^n(V)$
is one-dimensional and a non-zero vector is $x_1\wedge x_2\ldots\wedge x_n$. If $g\in G$ has the matrix $A$ with
respect to this basis, then
$$ g\cdot (x_1\wedge x_2\ldots \wedge x_n)=\det A\cdot x_1\wedge x_2\ldots \wedge x_n.$$

Let now $V$ be a supervector space of dimension $(m|n)$. The (super)determinant of 
an endomorphism of $V$ was introduced by Berezin. Fix a homogeneous basis of $V$,
$x_1,x_2\ldots x_{m+n}$ where the first $m$ elements are even and the rest are odd (such
a basis is called distinguished). Then an endomorphism of $V$, with respect to this basis,
has the matrix of the following form
$$Z=\left(\begin{array}{cc}A&B\\ C&D\end{array}\right)$$
where $A,D$ are square matrices of dimension $m\times m$ and $n\times n$, respectively, whose
entries are even, and $B,D$ are matrices of type $m\times n$ and $n\times m$, whose entries
are odd elements. Assume that $D$ is invertible in the usual sences, define the super determinant
of $Z$ to be
$$\sdet Z=\det T^{-1}\det(A-CD^{-1}B).$$
It is shown that the the matrix $Z$ is invertible iff its super determinant is and that the super determinant
is multiplicative. It is however not clear why the definition of $\sdet$ is independent of the choice of bases
(our basis is a distinguished basis).

In \cite{manin3} Manin suggested the following construction to define the super determinant. Let $V^*$
denote the  vector space to dual $V$ with the dual basis $\xi^1,\xi^2,\ldots,\xi^n$, $\xi^i(x_j)=\delta^i_j$.
Manin introduced the following Koszul complex; its $(k,l)$-term is given by $K^{k,l}:=\Lambda^k\ot S^{l*}$,
where $\Lambda^n$ and $S^n$ are the $n$-th homogeneous components of the exterior
and the symmetric tensor algebra over $V$. The differential $d^{k,l}:K^{k,l}\lora K^{k+1,l+1}$
is given by
$$d^{k,l}(h\ot \phi)=\sum_ihx_i\ot \xi^i\wedge\phi.$$
It is easy to check that $d_{k,l}$ is $G$-equivariant hence the homology groups of this complex
are representations of $G$. On the other hand, one can show that  this complex 
is exact every where except at the term $(m,n)$, where the homology group is one-dimensional, 
thus, it defines a one-dimensional representation of $G$. It turns out that elements of $G$ acts on
this representation by means of its super determinant, in other words, the definition of
 the super determinant  is basis free.

The quantum semigroup of type $A$ is the ``spectrum" of the quadractic algebra 
$$E:=k\langle z^i_j\rangle/(R^{ij}_{uv}z^u_kz^v_l=z^i_tz^j_sR^{ts}_{kl})$$
where $R$ is a Hecke symmetry (see \S\ref{sect1}). The Hecke symmetry resembles the usual flipping
operator $a\ot b\loma b\ot a$ or $a\ot b\loma (-1)^{\hat a\hat b}b\ot a$ ($a,b$ are homogeneous)
in super symmetry. 

In \cite{gur1,ls}, a Koszul complex is defined for $R$. For that, one first has to
define the quantum exterior and quantum symmetric tensor algebra by means of certain projectors
on $V^{\ot n}$. It is still an open question, whether this complex has the homology group concentrated
at a certain term and its dimension is one. Some efforts have been made.  Gurevich \cite{gur1}
showed this for even Hecke symmetries (i.e., those which induce finite-dimensional exterior algebra),
Lyubashenko and Sudbery \cite{ls} showed this for  Hecke sums of an odd and an even
Hecke symmetries.

In this paper, assuming that $R$ depends algebraically on $q$, where $q$ runs in $\mathbb C$,
  we give the affirmative answer to this question for an algebraically dense set of value of $q$.
Our tactic is first  to use a new result of  Deligne \cite{deligne02} to check the case $q=1$. Then 
using standard argument we show that for a dense set of values of $q$, the homology group of $K$
has the dimension less than that  of the  corresponding homology groups when $q=1$.
In other words, for an algebraically dense subset of  $\mathbb C$, the homology has dimension
at most 1. It remains to show the non-vanishing of the homology.
Taking the tensor product of the complex with a suitably chosen comodule we obtain a new complex
whose terms are all $E$-comodules. 
We decompose these comodule using the Littlewood-Richardson formula
and derive the non-vanishing of the homology.

\section{Hecke symmetries and the associated quantum groups}\label{sect1}
We work over an algebraically closed field $\k $ of 
characteristic zero. Let $V$ be a vector space over $\k $
 of finite dimension $d$. Let $R:V\ot V\lora V\ot V$ be an invertible 
operator. $R$ is called a {\it Hecke symmetry} if the following conditions are fulfilled:
\begin{itemize}
\item $R_1R_2R_1=R_2R_1R_2$, where $R_1:=R\ot \id_V, R_2:=\id_V\ot R$,
\item $(R+1)(R-q)=0$ for some $q\in\k $.
\item The half adjoint to $R$, $R^\sharp:V^*\ot V\lora V\ot V^*$;
 $\langle R^\sharp(\xi\ot v,w\rangle=\langle\xi,Rv\ot w\rangle$,
is invertible.\end{itemize}
Through out this work we will  assume that $q$ {\em is not a root of unity other then the unity itself}.
If $q=1$, $R$ is called vector symmetry. Vector symmetries were introduced by Lyubashenko \cite{lyu1}
and generalized to Hecke symmetries by Gurevich \cite{gur1}.

Let us fix a basis $x_1,x_2,\ldots,x_d$ of $V$. Then $R$ can be given in terms of a matrix,
also denoted by $R$, $R(x_i\ot x_j)=x_k\ot x_lR^{kl}_{ij}$, we adopt
the convention of summing up after the indices that appear both in the lower and upper
places. The matrix ${R^\sharp}_{ij}^{kl}$ is given by ${R^\sharp}^{kl}_{ij}=R_{jl}^{ik}$.
Therefore, the invertibility of $R^\sharp$ can be expressed as follows:
there exists a matrix $P$ such that 
$P^{im}_{jn}R^{nk}_{ml}=\delta^i_l\delta^k_j$.

To a Hecke symmetry $R$, there associated the following quadratic algebras:
\bbas S&:=&k\langle x_1,x_2,\ldots,x_d\rangle/(x_kx_lR^{kl}_{ij}=qx_ix_j)\\
\Lambda&:=&k\langle x_1,x_2,\ldots,x_d\rangle/(x_kx_lR^{kl}_{ij}=qx_ix_j)\\
E&:=&k\langle z^1_1,z^1_2,\ldots,z^d_d\rangle/(z^{i}_{m}z^j_nR^{mn}_{kl}=R^{ij}_{pq}z^p_kz^q_l)
\eeas
$S$ and $\Lambda$ are called respectively the function algebra and the exterior-algebra
on the corresponding quantum space and $E$ is called the function algebra on the
corresponding quantum endomorphism space or the matrix quantum semi-group.
 
  $E_R$ is a coquasitriangular bialgebra \cite{l-t,lyu}.
The coproduct on $E_R$ is given by $\Delta(Z)=Z\ot Z$. 
The coquasitriangular structure is given by $r(z^i_j,z^k_l)=R^{ki}_{jl}$.
$S$, $\Lambda$ and all their homonogeneous components are
right $E$-comodules. In particular, $V$ is a comodule and
the induced braiding on $V\ot V$ is exactly $R$.

$E_R$ is naturally $\mathbb{N}$-graded, let $E_n$ be its n-th homogeneous component. 
Then $E_n$ is a coalgebra and $V^{\ot n}$ is its comodule, hence an $E^*_n$-module.

On the other hand, $R$ induces a representation of the Hecke algebra $\H_n$ 
(see, e.g., \cite{dj1,dj2}) on $V^{\ot n}$,
 denoted by $\rho_n:\H_n\lora \End(V^{\ot n})$, for any $n\geq 2$. Explicitly, $\rho$ maps the
generator $T_i$ of $\H_n$ to the operator $R_i:=\id^{\ot i-1}\ot R\ot\id^{n-i-1}$.
We have the following ``Double centralizer theorem" \cite{ph97}
\begin{fact} 
The algebras $\rho_n(\H_n)$ and $E_n^*$ are centralizers of 
each other in $\End_{\k }(V^{\ot n})$.\end{fact}

Let $x_n\in\H_n$ be the central idempotent that induces the trivial representation,
$x_n=\sum_{w\in \Ss_n}(-q)^{-l(w)}T_w/[n]_{q}$, where $T_w$ are the generators of
$\H_n$ as a $k$-vector space, indexed by elements of the symmetric group $\Ss_n$.
 The operator $X_n:=\rho_n(x_n)$ is called the $q$-symmetrizer, 
it is a projection on $V^{\ot n}$. The projection $V^{\ot n}\lora S^n$ restricted
to $\Im X_n$ is an isomorphism. Since $R$ is a morphism of $E$-comodules
the above isomorphisms are isomorphism of $E$-comodules, too.

Analogously, let $y_n\in \H_n$ be the central idempotent that 
induces the signature representation of $\H_n$:
$y_n:=([n]_{1/q})^{-1}\sum_{w\in \Ss_n}(-q)^{-l(w)}T_w$. 
The operator $Y_n:=\rho_n(y_n)$ is called the $q$-anti-symmetrizer, 
it is a projection on $V^{\ot n}$. The projection $V^{\ot n}\lora \Lambda^n$ restricted 
to $\Im Y_n$ is an isomorphism.

A Hecke symmetry $R$ is called {\em even} (resp. {\em odd}) Hecke symmetry of rank $r$ 
iff $Y_{r+1}=0$, and $Y_r\neq 0$ (resp. $X_{r+1}=0$ and $X_r\neq 0$).
One can show that, in this case, $Y_r$ (resp. $X_r$) has rank 1 \cite{gur1}.

By means of the above double centralizer theorem, $E$ is cosemisimple 
(i.e. its comodules are semisimple) and simple $E$-comodules can be 
described by primitive idempotents of the Hecke algebras, thus, partitions.
For partitions $\lambda,\mu$ and 
$\gamma$, the multiplicity of $M_\gamma$ (the simple $E_R$
 comodule corresponding to $\gamma$) in $M_\lambda\ot M_\mu$
 does not depend on $R$, in fact, it is equal to the corresponding 
Littlewood-Richardson coefficient $c_{\lambda\mu}^\gamma=(s_\lambda s_\mu,s_\gamma)$, 
where $s_\lambda$ are the Schur functions (cf. \cite{mcd2}). Note that however that
not any partition defines a simple comodules, some of them may give zero-modules.

To have more precise information on the simple comodules of  $E$, we need
the Poincar\'e series of $S$, $\Lambda$. Using theory 
of symmetric functions \cite{mcd2} and a theorem of 
Edrei \cite{edrei} we have the following \cite{ph97c}:
\begin{fact}\label{s5} 
$P_\Lambda(t)$
 is a rational function having negative roots and positive poles.
Assume that  $P_\Lambda$ has $m$ roots and $n$ poles, then
simple $E_R$-comodules are parameterized by partitions 
$\lambda$ for which $\lambda_{m+1}\leq n$. 
 \end{fact}

\noindent{\sc Definition}. The pair $(m,n)$ is called
the birank of $R$.


\smallskip

\noindent {\sc Examples.} The following are so far examples of Hecke symmetries.
 \begin{itemize}
\item The solutions of the Yang-Baxter 
equation of series $A$, due to Drinfel'd and Jimbo \cite{jimbo86} is an example of even Hecke symmetries.
 The associated quantum groups are  called standard deformations of GL$(n)$.
\item Cramer and Gevais \cite{c-g} found another series of solution which are also even Hecke 
symmetries. 
\item Hecke sums of odd and even Hecke symmetries are 
examples of non-even, non-odd Hecke symmetries \cite{majid-markl}.
\item Takeuchi and Tambara found a Hecke symmetry which is neither even nor a Hecke sum of an 
odd and an even Hecke symmetries \cite{tt}.
\item Even Hecke symmetries of rank 2 was classified by Gurevich \cite{gur1}. He also show that on 
each vector space of dimension $\geq 2$, there exists an even Hecke symmetries of rank 2.
\item Hecke symmetry of birank $(1,1)$ was classified by the author \cite{ph99}.
\end{itemize}

\smallskip

The quantum group of type $A$ is define to be the ``spectrum" of the subsequently defined
Hopf algebra.
Let  $T=(t_i^j)$ be a $d\times d$ matrix of new variables. The
 Hopf algebra associated to $R$ is a factor algebra of the free associative
algebra over entries of $Z$ and $T$:
\begin{equation}\label{H_R}
 H_R:= T\langle Z,T\rangle \left/\left( RZ_1Z_2=Z_1Z_2R, TZ=ZT=\id            \right)\right.
\end{equation}
$H_R$ is a Hopf algebra, the antipode is given by $S(Z)=T$. 
The coquasitriangular struture on $E_R$ can be extended on to $H_R$
thanks the closedness of $R$. 

The structure of $H_R$-comodules is, in general, much more 
complicated then the one of $E_R$-comodules. The best handled 
case is when $R$ is an even Hecke symmetry, i.e., when $P_\Lambda(t)$
 is a polynomial. 
We have, however the following result \cite{ph97b}.
\begin{fact}\label{s18} The natural map $E_R\lora H_R$ is injective. 
Consequently, every simple $E_R$-comodule is a simple $H_R$-comodule. 
\end{fact}

Among $H_R$-comodules which are not $E$-comodules, the super determinant plays 
an important role. The well-known tool for defining the quantum super determinant serves 
the Koszul complex (of second type) introduced by Manin \cite{manin3}. 
This is a (bi-)complex, whose $(k,l)$ term is $\Lambda^k\ot S^{*l}$. 
The differential is induced from the dual basis map. The homology 
group of this complex is an $H_R$-comodule, if it is one dimensional 
over $k$, it defines a group-like element in $H_R$ called homological 
determinant or quantum super determinant or, in some cases, quantum Berezinian.

\section{The Koszul complex}
We begin with the description of the Koszul complex.
For convennience, we first  fix the following notion of the dual comodule 
of a tensor product of two or more comodules. 
For two (rigid) comodules $V,W$, the dual to $V\ot W$ is defined to be 
$W^*\ot V^*$ with the evaluation map $\ev_{V\ot W}=\ev_W\circ(W^*\ot \ev_V\ot W)$.
Analogously, one defines the dual to longer tensor products.

Fix a basis $x_1,x_2,\ldots,x_n$ of $V$ and let $\xi^1,\xi^2,\ldots,\xi^n$ be the dual
basis in $V^*$, we define the dual basis map $\db:k\lora V\ot V^*$, $\db(1)=\sum_ix_i\ot \xi^i$.
This map does not depend on the choice of basis.
The term $K^{k,l}$ of the Koszul complex associated to $R$ is
$\Lambda^k\ot S^{l*}$, the differential $d_{k,l}$  is given by:
\bbs \Lambda^k\ot S^{l*}\rightarrow V^{\ot l}\ot V^{*\ot l}\stackrel{\id\ot\db_V\ot \id}{\lora} 
V^{\ot k+1}\ot V^{*\ot l+1}\stackrel{Y_{k+1}\ot X_{l+1}^*}{\longrightarrow}
\Lambda^{k+1}\ot S^{l+1*},\ees
where $X_l, Y_k$ are the $q$-symmetrizer operators introduced in the
 previous section.
One defines another differential $d'$ as follows:
\bbs \Lambda^k\ot S^{l*}\rightarrow V^{\ot l}\ot V^{*\ot l}\stackrel{\id
\ot\ev_V\tau_{V\ot V^*}\ot\id}{\lora} V^{\ot k-1}\ot V^{*\ot l-1}
\stackrel{Y_{k-1}\ot X_{l-1}^*}{\longrightarrow}\Lambda^{k-1}\ot S^{l-1*},\ees
where $\tau_{V,V^*}$ denotes the braiding on $V\ot V^*$ induced from the coquasitriangular
structure on $H_R$, its matrix is given by $P$:$R_{im}^{jn}P_{nk}^{ml}=\delta_i^l\delta_k^j.$
Then $d$ and $d'$ satisfy \cite{gur1}
\bbs \left(qdd'+d'd\right)\left|_{K^{k,l}}\right.=q^k(\rk_qR+[l-k]_q)\id.\ees
where $\rk_qR:=P^{ij}_{ij}$.
Hence, if $\rk_qR\neq -[l-k]_q$, the cohomology group at the term $(k,l)$ vanishes. 
Thus, all complexes except at most one are acyclic.

\begin{thm} \label{non-vanishing}Let $R$ be a Hecke symmetry of birank $(m,n)$. Then $\rk_qR=-[n-m]_q$
and the homology of the Koszul complex at the term $(m,n)$ is non-vanishing.
Consequently, on the Hopf algebra $H$, there exists a non-zero integral and the simple
$H$-comodule $M_\lambda$ is injective and projective if and only if $\lambda_m\geq n$.\end{thm}
\proof
Since $R$ has birank $(m,n)$, simple $E$-comodules are parameterized
by partions satisfying $\lambda_{m+1}\leq n$. Using this fact and the Littlewood-Richardson
formula, we can easily show that
\bbas&& \Hom^E(M_{((n+1)^m)}\ot S^n,M_{(n^{m+1})}\ot \Lambda^m)=\k\\
&&\Hom^E(M_{((n+1)^m)}\ot S^{n-1},M_{(n^{m+1})}\ot\Lambda^{m+1})=0\\
&&\Hom^E(M_{((n+1)^m)}\ot S^{n+1},M_{(n^{m+1})}\ot\Lambda^{m+1})=0.\eeas
As a consequence, $M_{(n^{m+1)})}\ot\Lambda^m\ot S^{n*}$ contains
 $M_{((n+1)^m)}$ while the comodules
$M_{(n^{m+1})}\ot\Lambda^{m-1}\ot S^{n*}$, $M_{(n^{m+1})}\ot\Lambda^{m+1}
\ot S^{m+1*}$  do not. 

Assume that $\rk_qR\neq -[n-m]_q$. Then the complex is exact at $K^{m,n}$ 
and $dd'+d'd=q^m(\rk_qR+[n-m]_q)\id\neq 0$.

On the other hand, since $M_{(n+1)^m}$ cannot be a submodule of 
$M_{(n^{m+1})}\ot \Lambda^{m+1}\ot S^{n+1*}$, 
the restriction of $\id_{M_{((n+1)^m)}}\ot d^{m,n}$ on it should be zero. Analogously, the restriction
of $\id_{M_{((n+1)^m)}}\ot d^{m,n}$ on $M_{((n+1)^m)}$ is 0.
 Thus, the restriction of $dd'+d'd$ on $M_{(n+1)^m}$ 
is zero, a contradiction. Therefore $\rk_qR= -[n-m]_q$. 

 According to a result of \cite{ph98b} if  $\rk_qR=-[n-m]_q$ then $H$ possesses a non-zero  integral and
in this case, according to a result of \cite{ph99}, $M_\lambda$ is a splitting comodule (i.e.
injective and projetive in $H$-comod) iff $\lambda_m\geq n$. Thus, $M_{((n+1)^m)}$ is projective hence 
cannot be a subquotient of $M_{(n^{m+1})}\ot\Lambda^{m-1}\ot S^{n-1*}$, in particular, it cannot be a 
subcomodule of $M_{(n^{m+1})}\ot\Im d^{m-1,n-1}$. Therefore
$$M_{(n^{m+1})}\ot\Im d^{m-1,n-1}\neq M_{(n^{m+1})}\ot\Ker d^{m,n}.$$
Thus,  the sequence
$$M_{(n^{m+1})}\ot\Lambda^{m-1}\ot S^{n*}\rightarrow M_{(n^{r+1)})}
\ot\Lambda^m\ot S^{n*}  \rightarrow M_{(n^{m+1})}\ot\Lambda^{m+1}\ot S^{m+1*}$$
which is obtained  by tensoring $K^{\cdot\cdot}$ with $M_{(n^{m+1})}$
 is not exact at the term $(m,n)$, whence neither is $K^{\cdot\cdot}$.\eee

\section{The case $q=1$} Assume in this section, $q=1$, thus, $R^2=1$ and $H\comod$
is a tensor category (i.e., symmetric rigid monoidal). By a theorem of Deligne,
 there exists a faithful
and exact, tensor (i.e. symmetric monoidal) functor $\fun F$  from $H\comod$ to the category of 
vector superspaces. Under this functor, $V$ is mapped to a certain vector super space $\bar V$
and $R$ is mapped to the supersymmetry on $\bar V\ot \bar V$, denoted by $T$.

 We can therefore reconstruct a super bialgebra $\bar E$ and a 
Hopf super algebra $\bar H$ from $V$ and $T$.
We will show that this Hopf superalgebra is isomorphic to the function algebra over the
general linear supergroups $GL(n)$, where $(m,n)$ is the birank of $R$, or, in
other words, the super dimension of $\bar V$ is $(m|n)$. Indeed, $\bar E$ is the function algebra
on $\End(\bar V)$ and the image of  $M_\lambda$ under the embedding $\fun F$ are simple
$\bar E$-comodules. Since $\fun F$ is faithful and exact and since 
$M_\lambda\neq 0\Leftrightarrow \lambda_{m+1}< n$, 
we conclude that $\bar E$ is isomorphic to the function algebra 
on $M(m|n)$. Hence $\bar H$ is isomorphic to the function algebra on $GL(m|n)$, by virtue
of \ref{s18}.

Let ${\bar K}^{\cdot\cdot}$ denote the image of the complex ${K}^{\cdot\cdot}$.
Then the homology of ${\bar K}^{\cdot\cdot}$  is concentrated at the term $(m,n)$,
and is one-dimensional; it defines the super determinant. As a consequence, the homology
of ${ K}^{\cdot\cdot}$ is also concentrated at the term $(m,n)$, for $\fun F$ is faithful
and exact. Let $D$ denote the homology of ${ K}^{\cdot\cdot}$. Then $\bar D$, the image
of $D$ under $\fun F$, is one-dimensional and hence invertible, consequently,
$$\fun F(D^*\ot D)\cong \fun F(D^*)\ot \fun F(D)\cong {\bar D}^*\ot\bar D\cong k,$$
where the last isomorphism is given by the evaluation morphism, that is the image of
$\ev_D$ under $\fun F$.
Since $\fun F$ is faithful and exact, we conclude that $D^*\ot D\cong k$, that is
$D$ is invertible, hence one-dimensional.
$\bar M_\lambda$ denote the image
where $\lambda$ runs in the set of partitions for which
\begin{thm}\label{caseq=1} Let $R$ be a vector symmetry of birank $(m,n)$.
 Then the associated Koszul complex
is exact every where except at the term $(m,n)$ where it has a one-dimensional
homology group which determines a group-like element called homological
determinant.\end{thm}
\section{The case $q$ generic}
Using the result of the previous section we show in this section that given a Hecke
symmetry of birank $(m,n)$ that depends algebraically on $q$, then,
 for a dense set of values $q$, the associated Koszul complex is exact every where 
except at the term $(m,n)$, where is has a one-dimensional homology group and
thus determines a group-like element in $H$, called the homological determinant.
In this section $k$ will be assumed to be the field $\mathbb C$ of complex numbers.

Thus let $R=R_q$ be a Hecke symmetry depending on a parameter $q\in \mathbb C$.
We first observe that the dimension of $\Lambda_q^k$ does not depend on $q$, as far
as $q$ is not a root of unity. Indeed, $\Lambda_q^k$ can be defined as the
image of a projection, its dimension can be given as the trace of a matrix which
depend algebraically on $q$, since $\mathbb C$ substracted the set of root
of unity is still connected, we conclude that this trace, being always integral must be a
constant. The same happens with $S_q^l$.
Thus, the terms of $K^{\cdot\cdot}$ has the dimension not depending on $q$.

On the other hand, observe that the rank of the operator $d^{k,l}_q$, for almost any $q$ (that 
is except a finite number of values of $q$)  is large
then the rank of $d^{k,l}_1$ and for the kernel of $d_q^{k,l}$ we have the reversed
inequality. Consequently, the dimension over $\k$ of the homology group
$H(K_q^{k,l})$ for almost  any $q$ is lest then or equal to the dimension of 
$H(K_q^{k,l})$. According to Theorems \ref{non-vanishing} and \ref{caseq=1}, 
we conclude that for an algebraically dense set
of values of $q$, $H(K_q^{k,l})=0$, for all $(k,l)\neq (m,n)$ and 
$H(K_q^{m,n})=k$.
\begin{thm} Let $R=R_q$ be a Hecke symmetry over $\mathbb C$, depending
algebraically on $q$. Then there is an algebraically dense set of values of $q$
for which the homology of the Koszul complex is one-dimensional and concentrated
at the term $(m,n)$, where $(m,n)$ is the birank of $R$.\end{thm}

\end{document}